\title{\vspace{-0.6cm} On graphs decomposable into induced matchings of linear sizes}

\documentclass[11pt]{article}

\usepackage{amsmath, amssymb, amsthm, url}

\newenvironment{pf}
     {\medskip\noindent{\bf Proof.}\hspace{1mm}}
      {\hfill$\Box$\medskip}

\newtheorem{theorem}{Theorem}[section]
\newtheorem{proposition}[theorem]{Proposition}
\newtheorem{lemma}[theorem]{Lemma}

\newtheorem{corollary}[theorem]{Corollary}

\newtheorem{problem}[theorem]{Problem}

\newcommand{\dist}{{\textup{dist}}}

\date{}

\author{
Jacob Fox \thanks{
Department of Mathematics, Stanford University, CA 94305-2125.
Email: jacobfox@stanford.edu. Research supported
by a Packard Fellowship, by NSF Career Award DMS-135212, and by an Alfred P. Sloan Fellowship.
}
\and
Hao Huang \thanks{
Department of Math and CS, Emory University, Atlanta, GA 30322.
Email: hao.huang@emory.edu.
}
\and
Benny Sudakov \thanks{Department of Mathematics, ETH, 8092 Zurich, Switzerland. Email: benjamin.sudakov@math.ethz.ch. Research supported in part by SNSF grant 200021-149111.}
}

\newenvironment{rem}{\noindent\textbf{Remark.}}{\medskip}

\begin{document}
\maketitle
\abstract
We call a graph $G$ an $(r,t)$-Ruzsa-Szemer\'edi graph if its edge 
set can be partitioned into $t$ edge-disjoint induced matchings, each of size $r$. These graphs were introduced in 1978 and has been extensively studied since then. In this paper, we consider the case when $r=cn$. For $c>1/4$, we determine the maximum possible $t$ which is a constant depending only on $c$. On the other hand, when $c=1/4$, there could be as many as $\Omega(\log n)$ induced matchings. We prove that this bound is tight up to a constant factor. Finally, when $c$ is fixed strictly between $1/5$ and $1/4$, we give a short proof that the number $t$ of induced matchings is $O(n/\log n)$. We are also able to further improve the upper bound to $o(n/\log n)$ for fixed $c> 1/4-b$ for some positive constant $b$.


\section{Introduction} \label{section_introduction}
We call a graph $G=(V, E)$ an $(r, t)$-Ruzsa-Szemer\'edi graph, or $(r, t)$-RS graph for short, if its edge set $E(G)$ can be partitioned into $t$ pairwise edge-disjoint
induced matchings $M_1, \cdots, M_t$, each consists of $r$ edges.
These graphs were first introduced in the famous paper by Ruzsa and Szemer\'edi \cite{ruzsa-szemeredi}, in which they show that there is no $n$-vertex $(r,t)$-RS graph for $r, t$ both linear in $n$. 
They also proved that this result implies the celebrated theorem of Roth \cite{roth}, that a subset $S$ of $[n]=\{1, \cdots,n\}$ without nontrivial $3$-term arithmetic progressions has size at most $o(n)$. Ruzsa and Szemer\'edi also constructed $(n/e^{O(\sqrt{\log n})}, n/3)$-RS graphs on $n$ vertices, based on a result of Behrend \cite{behrend} that there is a surprisingly large subset of $[n]$ without non-trivial $3$-term arithmetic progressions. They used these graphs together with the regularity lemma of Szemer\'edi \cite{regularity} to solve an extremal problem of Brown, Erd\H os, and S\'os \cite{brown-erdos-sos-1, brown-erdos-sos-2}, showing that the maximum number of edges in a $3$-uniform $n$-vertex hypergraph in which no $6$ vertices span at least $3$ edges is at least $n^{2-o(1)}$ and at most $o(n^2)$ for sufficiently large $n$. In this paper we consider the following natural question.

\begin{problem}
\label{RS-problem}
For which values of $r$ and $t$ does there exist an $(r,t)$-RS graph on $n$ vertices?
\end{problem}

This question has been studied in depth for decades, and has found several applications in combinatorics, complexity theory and information theory (see e.g. \cite{alon, alon-moitra-sudakov, alon-shapira, alon-shapira-2, hastad-wigderson, fischer-et-al}). One line of research was to find very dense RS-graphs decomposable into large induced matchings. 
An early result of this type by Frankl and F\"uredi \cite{frankl-furedi} implies that for any fixed $r$ there are $(r, t)$-RS graphs on $n$ vertices with $rt=(1-o(1)) \binom{n}{2}$ edges. However the techniques employed in their proof cannot provide induced matchings of size larger than $\Theta(\log n)$. Later Birk, Linial and Meshulam \cite{birk-linial-meshulam} noticed that such very dense $(r,t)$-RS graphs with $r \rightarrow \infty$ are useful for designing a communication protocol over a shared directional multi-channel. For this application, they construct $(r,t)$-RS graphs on $n$ vertices with $r=(\log n)^{\Omega(\log \log n/(\log \log \log n)^2)}$ and $rt$ is roughly equal to $n^2/24$. Note that none of the constructions mentioned so far gives an $n$-vertex $(r, t)$-RS graph with positive edge density and at the same time being an edge-disjoint union of induced matchings of size $n^\varepsilon$ for some constant $\varepsilon$. This range of parameters is important for some applications when there is a tradeoff between the number of missing
edges and the number of induced matchings. Indeed Meshulam \cite{meshulam} conjectures that there exist no such
graphs. This conjecture was recently disproved by Alon, Moitra and Sudakov \cite{alon-moitra-sudakov} with a surprising construction of graphs with
edge density $1-o(1)$ that can be partitioned into induced matchings of size $n^{1-o(1)}$, which is nearly linear in $n$. Their constructions have also found a couple of interesting applications for communication problems over shared channels, linearity testing, communication complexity and the directed Steiner tree problem.

Another appealing direction is to determine the maximum number of induced matchings when their sizes are linear in $n$. 
Fischer et al. in \cite{fischer-et-al} gave a family of $(n/6-o(n), n^{\Omega(1/\log \log n)})$-RS graphs. They use these graphs to establish
an $n^{\Omega(1/\log \log n)}$ lower bound for testing monotonicity in general posets. A slight twist of their construction gives
$(cn, n^{\Omega(1/\log \log n)})$-RS graphs for any positive constant $c < 1/4$. In this paper, we first study the range $c > 1/4$.
It was noticed earlier by Alon \cite{alon_pc}  that for this range, there can only be constantly many matchings. Our first result determines the maximum number of induced matchings for $c>1/4$. Its proof combines techniques from coding theory together with properties of Kneser graphs.

\begin{theorem}\label{double_count}
Suppose $G$ is an $(r,t)$-RS graph on $n$ vertices, then
\begin{equation} \label{eqn1}
r \leq
\begin{cases}
\dfrac{n}{4}\left(1+\dfrac{1}{t}\right) & \text{if $2 \nmid t$,}\\
\dfrac{n}{4}\left(1+\dfrac{1}{t+1}\right)  & \text{if $2 \mid t$}
\end{cases}. 
\end{equation}
Moreover, these bounds are tight for every positive integer $t$ and infinitely many $n$.
\end{theorem}
  
For $r=n/4$, it was known that there are RS-graphs with logarithmically many induced matchings \cite{rashod, alon_pc} but there were no good upper bounds in this case. Our second theorem shows that up to a constant factor the maximum number of induced matchings when $r=n/4$  
is indeed logarithmic.

\begin{theorem} \label{main_theorem}
If an $n$-vertex graph $G$ is an $(n/4, t)$-RS graph, then $t \le (6+o(1)) \log_2 n $.
\end{theorem}

The Problem \ref{RS-problem} becomes much more difficult when the size $r=cn$ of the induced matching is below the threshold $n/4$.
For this range, the best known upper bound
is due to Fox \cite{fox}, who improved estimates which one obtaines from the Szemer\'edi regularity lemma \cite{regularity}. His result implies that the number of matchings is at most $n/\log^{(x)}n=o(n)$, where $x=O(\log(1/c))$ and $\log^{(x)} n$ denotes the $x$-fold iterated logarithm. This shows that the number of induced matchings has to be sublinear.
On the other hand, the example of Fischer et al. \cite{fischer-et-al} shows that for $c<1/4$, the number of matchings could be $n^{\Omega(1/\log \log n)}$, which is much larger than $\log n$ but smaller than $n^{\varepsilon}$.  It would be very interesting to close the huge gap between the upper and lower bounds for some ranges of $c<1/4$. Here we are able to make some progress on this problem and to improve the upper bound to $O(n /\log n)$, when $c$ is fixed strictly between $1/5$ and $1/4$. 

\begin{theorem} \label{theorem_below_n/4}
For every $\varepsilon>0$, if $G$ is an $(r,t)$-RS graph on $n$ vertices with $r=cn$ for $1/5+\varepsilon \le c <1/4$, then $t=O(n/\log n)$. 
\end{theorem}

By a more careful analysis, we are able to further improve the upper bound to $o(n/\log n)$, when $c$ is slightly below $1/4$. To be more precise, the little-o term comes from the best known bound in the triangle removal lemma proved in \cite{fox}.
 
\begin{theorem} \label{theorem_just_below_n/4}
There exists an absolute constant $b>0$, such that for $r \ge (1/4-b)n$, if $G$ is an $(r,t)$-RS graph $G$ on $n$ vertices, then $t=n/((\log n)2^{\Omega(\log^* n)})=o(n \log n)$.
\end{theorem}

The rest of the paper is organized as follows. In the next section we describe a counting proof of the upper bound in Theorem \ref{double_count}, and verify its tightness by a construction based on Kneser graphs. In Section \ref{section_r=n/4}, we prove Theorem \ref{main_theorem}
and discuss potential improvements. Section \ref{section_r<n/4} consists of the proof of Theorem \ref{theorem_below_n/4} and \ref{theorem_just_below_n/4}, and some remarks on the case of regular graphs. The last section contains some concluding remarks and open problems.

\section{The size of matchings is $cn$ with $c>1/4$} \label{section_r>n/4}
In this section, we determine the maximum number of induced matchings in an $(r, t)$-RS graph $G$ for $r=cn$ and $c>1/4$. To prove an upper bound, we use ideas similar to the Plotkin bound \cite{plotkin} in coding theory, with a small twist. The lower bound construction is based on the properties of Kneser graphs.
\medskip

\noindent \textbf{Proof of Theorem \ref{double_count}:}
Suppose the edge set of $G$ can be partitioned into induced matchings $M_1, \cdots, M_t$, each containing exactly $r$ edges.
Denote by $V_i$ the set of vertices contained in the edges of $M_i$. Then $|V_i|=2r$.
Moreover, each of the $r$ edges of $M_i$ intersects $V_j$ in at most one vertex, 
since otherwise $V_i$ and $V_j$ must span a common edge of $G$. This implies that $|V_i \cap V_j| \le r$. Let $v_i \in \{0,1\}^n$ be the characteristic vector of $V_i$. Then
for all $1 \leq i< j \leq t$, the Hamming distance satisfies
$$\dist(v_i,v_j)=|V_i|+|V_j|-2|V_i \cap V_j| \ge 2r + 2r -2r = 2r.$$
This is already enough to show that $t$ is constant, using bounds from coding theory. But one can do slightly better. Let $v_0$ be the all-zero vector. To get a tight upper bound, notice that the above inequality can be extended to all $0 \leq i< j \leq t$ since, for $1 \le i \le t$, $|V_i|=2r$.
Denote by $a_i$ (resp. $b_i$) the number of vectors $v_j$ equal to $0$ (resp. $1$) in the $i$-th coordinate, then $a_i+b_i=t+1$. By double counting,
\begin{equation}\label{eqn2}
\begin{split}
2r \binom{t+1}{2} &\leq \sum_{0 \leq i <j \leq t} \dist(v_i, v_j)\\
&=\sum_{i=1}^n a_i b_i\\
&\leq \begin{cases}
n(t+1)^2/4 & \text{if $2 \nmid t$,}\\
nt(t+2)/4 & \text{if $2 \mid t$}
\end{cases}
\end{split}
\end{equation}
The last inequality follows from that $a_ib_i$ is maximized when $a_i=b_i=(t+1)/2$ for odd $t$, and $\{a_i,b_i\}=\{(t+2)/2, t/2\}$ for even $t$.
By simplifying the inequality we immediately obtain \eqref{eqn1}.

To show that the bound is tight, it suffices to consider the case $t=2k+1$. Let $H$ be $KG(2k+1, k)$, the Kneser graph whose vertices correspond to
all the $k$-subsets of a set of $2k+1$ elements, and where two vertices are adjacent if and only if the two corresponding sets are disjoint.
We define the matchings $M_1, \cdots, M_{2k+1}$ in the following way:
the edge $(A, B)$ belongs to $M_i$ if and only if $A \cap B= \emptyset$ and $A \cup B=[2k+1] \setminus \{i\}$.
It is easy to see that $B$ is determined after fixing $A$ and $i$,
which implies that $M_i$ forms a matching. In order to show that every matching is induced in $H$, we take $(A,B)$ and $(C,D)$ both from $M_i$ with $A \neq C, D$,
then $A \cup B = C \cup D = [2k+1] \setminus \{i\}$, it is not hard to check that $A \cap C$, $A \cap D$ are both nonempty and therefore $(A, C)$ and $(A, D)$ are not contained in any $M_j$. By calculation, $n = \binom{2k+1}{k}$, while $r=\frac{1}{2}\binom{2k}{k} = \frac{n}{4} (1+ \frac{1}{2k+1})$. Hence
$H$ is a $(\frac{n}{4}(1+\frac{1}{2k+1}), 2k+1)$-RS graph on $n$ vertices.
\qed

\medskip
\begin{rem}
By taking $m$ vertex-disjoint copies of the Kneser graph $KG(2k+1,k)$, we can construct an infinite family of $n$-vertex $(r,2k+1)$-RS graphs with $r=n(1+1/t)/4$. Moreover, solving $n=\binom{2k+1}{k}$ gives $k \sim \frac{1}{2}\log_2 n + \Theta(\log \log n)$, therefore we obtain
a $(\frac{n}{4}+\Theta(\frac{n}{\log n}), \log_2 n + \Theta(\log \log n))$-RS graph $H$ on $n$ vertices.
\end{rem}

\section{The size of matchings is $n/4$} \label{section_r=n/4}
From Theorem \ref{double_count}, we know that in a Ruzsa-Szemer\'edi graph, if every matching has size $r=cn$ where $c>1/4$, then the number of matchings is at most a constant which depends on $c$. On the other hand, when $c=1/4$, we saw in the end of the previous section that one can have logarithmic number of induced matching. The following very natural construction \cite{rashod, alon_pc} does better and gives $2 \log_2 n$ induced matchings. 

\begin{proposition}
There exists $(n/4, 2 \log_2 n)$-RS graph for every integer $n$ that is a power of $2$.
\end{proposition}
\begin{pf}
Let $k=\log_2 n$ and consider the $k$-dimensional hypercube graph $H$ with vertex set $\{0,1\}^k$, where two vectors are adjacent if their Hamming distance is $1$. We first partition the vertices into odd and even vectors, according to the parity of the sum of their coordinates.
For $1 \le i \le k$, we let the $i$-th matching $M_i$ consist of edges between vectors $\vec{v}$ and $\vec{v}+\vec{e}_i$, such that $\vec{v}$ is even and its $i$-th coordinate is $0$; and the $(k+i)$-th matching $M_{k+i}$ consist of edges between an odd vector $\vec{v}$ whose $i$-th coordinate equals $0$ and the vector $\vec{v}+\vec{e}_i$. This construction gives $2k=2\log_2 n$ matchings, and obviously each matching involves exactly half of the vertices. In order to verify that the matchings are induced, we consider two distinct edges from $M_i$, which are $(\vec{u}, \vec{u}+\vec{e}_i)$ and $(\vec{v}, \vec{v}+\vec{e}_i)$, such that both $\vec{u}$ and $\vec{v}$ are even and their $i$-th coordinates are $0$. Clearly the pairs $(\vec{u}, \vec{v})$ and $(\vec{u}+\vec{e}_i, \vec{v}+\vec{e}_i)$ cannot form edges of $H$ since they have the same parity. Moreover, $\vec{u}$ and $\vec{v}+\vec{e}_i$ (similarly, $\vec{v}$ and $\vec{u}+\vec{e}_i$) differ in at least two coordinates. Therefore for all $1 \le i \le k$, the matchings $M_i$ we defined are induced. Using a similar argument, we can also show that the matchings $M_{i+k}$ are induced.
\end{pf}

When $\log_2 n$ is an even integer, we can improve slightly the above construction by adding two additional induced matchings.
The first one consists of $(\vec{u}, \vec{v})$ where $\vec{u}+\vec{v}=\vec{1}$ and both are even vectors. The second matching contains all edges $(\vec{u}, \vec{v})$ where $\vec{u}+\vec{v}=\vec{1}$ and both are odd vectors. This gives the following corollary.
\begin{corollary}\label{sth}
There exists $(n/4, 2(\log_2 n +1))$-RS graphs on $n$ vertices for every $n$ that is an even power of $2$.
\end{corollary}

One should naturally ask how good are these constructions, i.e.,  what upper bound can we prove on the number 
$t$ of induced matchings in an $n$-vertex $(r,t)$-RS-graph with $r=\frac{n}{4}$. Note that in the proof of Theorem \ref{double_count}, when $r=n/4$, the inequalities in \eqref{eqn2} are always true for every $t$. Indeed, in this case there might be linearly many $n$-dimensional binary vectors such that their pairwise Hamming distance is equal to $n/2$ by modifying the well-known Hadamard matrix into a $(0, 1)$-matrix and taking the row vectors. In the rest of this section, we develop a different approach which shows that when every induced matching has size exactly $n/4$, there can be no more than $O(\log n)$ matchings. Therefore the above constructions are best possible up to a constant factor. Below is a quick outline of our proof.

First we show that in any graph $G$ which is decomposable into induced matchings of size $n/4$, there are many pairs of adjacent vertices whose sum of degrees is large. This gives an auxiliary graph $H$ which is the subgraph of $G$ whose edges have large degree sum. Using additional properties of induced matchings, we can show that $H$ has nice expansion properties. From there we can easily derive some estimates on the number of induced matchings. To illustrate this idea, we first prove a slightly weaker upper bound. The following lemma allows us to only consider bipartite graphs for this case.

\begin{lemma} \label{lemma_bipartite}
If $G$ is an $(r, t)$-RS graph on $n$ vertices, then its bipartite double cover $G \times K_2$ is an $(2r, t)$-RS graph on $2n$ vertices.
\end{lemma}
\begin{pf}
Denote by $G' = G \times K_2$ the bipartite double cover of the graph $G$. The vertices of $G'$ are $(v, i)$ with $v \in V(G)$ and $i \in \{0, 1\}$. Two vertices $(u, 0)$ and $(v, 1)$ are adjacent whenever $u$ and $v$ form an edge in $G$. Note that an induced matching $M_i=\{(u_j, v_j)\}_{j=1}^{r}$ in $G$  corresponds to a matching $M'=\{((u_j, 0), (v_j, 1))\}_{j=1}^{r} \cup \{((u_j, 1), (v_j, 0))\}_{j=1}^{r}$ in $G'$, which is of size $2r$.
It is also straightforward to check by definition that $M'$ is an induced matching. Therefore $G'$ is an $(2r, t)$-RS graph on $2n$ vertices,
\end{pf}

\begin{theorem} \label{main_theorem_weaker}
If $G$ is an $(\frac{n}{4}, t)$-RS graph on $n$ vertices, then $t \le 8 (\log_2 n+1) $.
\end{theorem}
\begin{pf}
From Lemma \ref{lemma_bipartite}, it suffices to show that for all $n$-vertex bipartite graphs $G$ whose edges can be decomposed into induced matchings $M_1, \cdots, M_t$, each of size $n/4$, $t$ is at most $8 \log_2 n$.

Denote by $d_v$ the degree of vertex $v$ in $G$. We consider the subgraph $H$ of $G$, with edges being the pair of vertices $u, v$ such that $d_u + d_v \ge t$ and $(u, v) \in E(G$). Note that $e(G)=\frac{nt}{4}$. By the Cauchy-Schwarz inequality,
\begin{align} \label{inequality_cs}
\sum_{(u, v) \in E(G)} (d_u + d_v - t) &= \left(\sum_{v \in V(G)} d_v ^2\right) - t\cdot e(G)
\ge n \left( \frac{\sum_{v \in V(G)} d_v}{n}\right)^2 - \frac{nt^2}{4} \nonumber\\
&=n \left(\frac{ nt/2}{n}\right)^2-\frac{nt^2}{4} =0.
\end{align}

For any edge $(u, v) \in E(G)$, all the edges incident to either $u$ or $v$ must belong to different matchings. Therefore $d_u + d_v \le t+1$.
If we denote by $E_i$ the number of edges $(u, v)$ such that $d_u + d_v = t+ i$, we have
$E_1+ E_0 + E_{-1} + \cdots + E_{-t} =nt/4$, while inequality \eqref{inequality_cs} implies that
$E_1 - \sum_{j=1}^{t} j E_{-j} \ge 0$. Summing these two inequalities gives $2E_1 + E_0 \ge nt/4$ and so $E_1 + E_0 \ge nt/8$. So $H$ is a $n$-vertex graph with at least $nt/8$ edges and thus its average degree is at least $t/4$. Hence, $H$ has a subgraph $F$ of minimum degree at least $t/8$.

Set $s=t/8$. For each vertex $v$ of $G$, let $A_v$ denote the set of induced matchings containing $v$. Clearly $|A_v|=d_v$.
We claim that if $v$ and $u$ are at distance $k$ in $F$, then when $k$ is odd, $|A_u \cap A_v| \le k$; and when $k$ is even, $|A_u \cap A_v^c| \le k$. This statement can be proved using induction. The base cases when $k=0$ and $1$ are obvious. Now we assume that it is true for all $k \le i$. For $k= i+1$, suppose $u$ and $v$ are at distance $k$. Let $w$  be a vertex at distance $1$ from $v$ and $i=k-1$ from $u$. When $i$ is odd, from the inductive hypothesis we have
$|A_v^c \cap A_w^c|=t-|A_v \cup A_w|=t-|A_v|-|A_w|+|A_v \cap A_w| \le 1$ and $|A_w \cap A_u| \le i$. Therefore,
$$|A_u \cap A_v^c| \le |A_u \cap A_w| + |A_w^c \cap A_v^c| \le i+1.$$
Similarly, when $i$ is even, we have $|A_v \cap A_w| \le 1$ and $|A_w^c \cap A_u| \le i$, and hence
$$|A_u \cap A_v| \le |A_v \cap A_w|+ |A_w^c \cap A_u| \le i+1.$$

Now choose an arbitrary vertex $v$ in $F$, the degree of $v$ in $F$ is at least $s=t/8$. For every integer $i \ge 0$, let $N_i$ be the set of vertices at distance $i$ from $v$ in graph $F$. By the assumption that $G$ is bipartite, each $N_i$ induces an independent set.
We denote by $e_{i}$ the number of edges of $F$ that are between $N_i$ and $N_{i+1}$ and contained in matchings in $A_v$ (resp. $A_v^c$) when $i$ is odd (resp. even). For odd $i$, we estimate the number of edges of $F$ between $N_i$ and $N_{i+1}$ that are contained in matchings in $A_v^c$ in two different ways.
Since every vertex $u \in N_i$ is contained in at least $s-|A_u \cap A_v| \ge s-i$ such edges, and every vertex $w \in N_{i+1}$ is contained in no more than
$|A_w \cap A_v^c| \le i+1$ such edges,
we have
$$(s-i) |N_i|-e_{i-1} \le (i+1) |N_{i+1}|-e_{i+1}.$$

Similarly when $n$ is even, by bounding the number of edges between $N_i$ and $N_{i+1}$ that belong to matchings in $A_v$, we obtain the same inequality as above.

Summing up the inequalities for $i=0, \cdots , k$, we have
$$\sum_{i=0}^k (s-i) |N_i|- \sum_{i=0}^{k-1} e_i \le \sum_{i=1}^{k+1} i |N_i| - \sum_{i=1}^{k+1} e_i.$$

Simplifying this inequality gives
$$(k+1)|N_{k+1}| \ge \sum_{i=0}^{k} (s-2i)|N_i|-e_0+e_{k+1}+e_k \ge \sum_{i=0}^{k} (s-2i)|N_i|.$$
The second inequality follows from the observation that $e_0=0$ since all edges between $N_0$ and $N_1$ are in $A_v$.

In the next step, we prove by induction that $|N_i| \ge \binom{s}{i}$. For $i=0$ and $1$ this is obvious. Now, assuming it is true for all $i \le k$, we have
\begin{align*}
|N_{k+1}| &\ge \frac{1}{k+1} \sum_{i=0}^{k} (s-2i)|N_i|
 \ge \frac{1}{k+1} \sum_{i=0}^{k} (s-2i) \binom{s}{i}\\
& = \frac{1}{k+1} \left( \sum_{i=0}^{k} (s-i) \binom{s}{i} - \sum_{i=0}^{k} i \binom{s}{i} \right)\\
& = \frac{1}{k+1} \left( \sum_{i=0}^{k} s \binom{s-1}{i} - \sum_{i=0}^{k} s \binom{s-1}{i-1} \right)\\
& =\frac{s}{k+1} \binom{s-1}{k} = \binom{s}{k+1}.
\end{align*}
Note that the number of vertices in $N_0 \cup N_1 \cup \cdots \cup N_s$ is at most $n$. We therefore have
$$n \ge \sum_{k=0}^{s} |N_i| \ge \sum_{k=0}^{s} \binom{s}{k} = 2^s,$$
Solving this inequality gives $s \le \log_2 n$ and hence $t \le 8 \log_2 n$.
\end{pf}

As we mentioned in the introduction, this bound can be further improved to $(6+o(1)) \log_2 n$ by the following modification of the above proof.\\

\noindent \textbf{Proof of Theorem \ref{main_theorem}:}
Recall that the proof above uses the fact that $H$ contains $E_1+E_0$ edges and therefore it contains a subgraph $F$ with minimum degree at least
$(E_1+E_0)/n$, and later on we obtained the inequality $n \ge 2^{(E_1+E_0)/n}$.
Let $H'$ be the subgraph of $G$ consisting of edges $(u, v)$ such that $d_u + d_v = t+1$, so $e(H')=E_1$. Graph $H'$ contains a subgraph $F'$
of minimum degree at least $E_1/n=s'$. By an argument similar to the one in the previous proof, we can show that if a vertex $u$ is at distance $k$ from $v$, then
$|A_u \cap A_v| \le (k+1)/2$ if $k$ is odd, and $|A_u \cap A_v^c| \le k/2$ if $k$ is even. Therefore, if we fix a vertex $v$ and let $N'_i$ be the set
of vertices having distance $i$ from $v$, and define $e'_i$ in a similar way, we have
$$(s'-\lceil i/2 \rceil)|N'_i|-e'_{i-1} \le \lceil (i+1)/2 \rceil |N'_{i+1}| -e'_{i+1}.$$
Summing up the inequalities from $i=0$ to $k$, and using $e_{-1}'=e_0'=0$ and $e_k',e_{k+1}' \geq 0$, gives
$$\lceil (k+1)/2 \rceil |N'_{k+1}| \ge \sum_{i=0}^k (s'-2\lceil i/2\rceil)|N'_i|.$$
We can use induction to show that $|N'_{2i}| \ge \binom{s'}{i}\binom{s'-1}{i}$ and $|N'_{2i+1}| \ge \binom{s'}{i+1}\binom{s'-1}{i}$. The $i=0$ case follows from $|N'_0| \ge 1$ and $|N'_1| \ge s'$. Moreover, using the inductive hypothesis, 
\begin{align*}
i \cdot |N'_{2i}| &\ge \sum_{j=0}^{i-1} (s'-2j)|N'_{2j}| + \sum_{j=0}^{i-1}(s'-2(j+1))|N'_{2j+1}|\\
&\ge \sum_{j=0}^{i-1} (s'-2j)\binom{s'}{j}\binom{s'-1}{j} + 
\sum_{j=0}^{i-1}(s'-2(j+1))\binom{s'}{j+1}\binom{s'-1}{j}
\end{align*} 
Substituting in the identity $(s'-j)\binom{s'}{j}\binom{s'-1}{j}=(j+1)\binom{s'}{j+1}\binom{s'-1}{j}$, and then the identity $(s'-j-1)\binom{s'-1}{j}=(j+1)\binom{s'-1}{j+1}$, and finally computing a telescoping sum, we obtain
\begin{align*}
i \cdot |N'_{2i}| &\ge \sum_{j=0}^{i-1} (s'-j-1)\binom{s'}{j+1}\binom{s'-1}{j}-j\binom{s'}{j}\binom{s'-1}{j}\\
&= \sum_{j=0}^{i-1} (j+1)\binom{s'}{j+1}\binom{s'-1}{j+1}-j\binom{s'}{j}\binom{s'-1}{j}\\
&=i \binom{s'}{i}\binom{s'-1}{i}
\end{align*} 
Similarly, it can be verified that $|N'_{2i+1}| \ge \binom{s'}{i+1}\binom{s'-1}{i}$ by the inductive hypothesis.

By counting the number of vertices in $F'$, we have
\begin{align*}
n &\ge  |N'_0|+\left(\sum_{i=1}^{s'-1} |N'_{2i-1}|+|N'_{2i}|\right) + |N_{2s'}| \\
&\ge 1+ \left(\sum_{i=1}^{s'-1} \binom{s'}{i}\binom{s'-1}{i} + \binom{s'}{i}\binom{s'-1}{i-1}\right) +1\\
&=1+ \left(\sum_{i=1}^{s'-1} \binom{s'}{i}^2\right) + 1 = \binom{2s'}{s'}.
\end{align*}
Therefore $$E_1/n = s' \le (\frac{1}{2}+o(1)) \log_2 n.$$
Together with the previously established inequalities $(E_1+E_0)/n = s \le \log_2 n$ and  $2E_1 + E_0 \ge nt/4$, we have
$$nt/4 \le 2E_1 + E_0 \le n(E_1/n + (E_1+E_0)/n) \le (\frac{3}{2}+o(1))n \log_2 n,$$
and $t \le (6+o(1)) \log_2 n$, which completes the proof. \qed

\medskip
\noindent \textbf{Remark.} If the $(n/4, t)$-RS graph $G$ is regular, then we can improve the upper bound to $t \le 2(\log_2 n + 1)$. In the proof of Theorem \ref{main_theorem_weaker}, the auxiliary graph $H$ is the same with $G$, in which every vertex has degree exactly $t/2$. Therefore for bipartite graph $G$, we have $n \ge 2^{t/2}$ and $t \le 2 \log_2 n$, and for general graphs we have $t \le 2(\log_2 n+1)$. This bound is the the best we could hope for by Corollary \ref{sth}.

\section{The size of matchings is below $n/4$} \label{section_r<n/4}
For the purpose of simplicity, we will use $\log$ in place of $\log_2$ throughout this section.\\

\noindent \textbf{Proof of Theorem \ref{theorem_below_n/4}:} Let $G=(V,E)$ and $E(G)$ be the disjoint union of $t$ induced matchings from $\mathcal{M}$ such that each matching consists of $r=cn$ edges with $1/5 + \varepsilon \le c < 1/4$. Suppose for contradiction $t=K n/\log n$ and $K=100/\varepsilon$. We start by defining $V_0=V(G)$. For each positive integer $i$, let $v_i$ be a vertex of maximum degree in the induced subgraph of $G$ on the vertex set $V_{i-1}$, $N_i$ be the set of neighbors of $v_i$ in $V_{i-1}$, and $V_i = V_{i-1} \setminus N_i$. Suppose $|N_i|=n_i$, then by definition $n_1 \ge n_2 \ge \cdots$ and $n_i \le d_i$, where $d_i$ is the degree of $v_i$ in $G$. Let $k$ be the maximum integer
such that $n_k \ge 2n/\log n$. Since all the $N_i$'s are disjoint, we have $n \ge \sum_{i=1}^k |N_i| \ge k \cdot (2n/\log n)$, hence $k \le (\log n)/2$.
Note that $|V_i|=n-\sum_{j=1}^i n_j$, and its number of edges $e(V_i)$ is at least $e(G)- \sum_{j=1}^i n_j^2$, since every vertex in $N_j$ has degree at most $n_j$ in $V_{j-1}$. By the choice of $k$, we also know that $e(V_k) \le 1/2 \cdot (2n/\log n) \cdot n = n^2 /\log n$.

Here is the key observation: for every $i$, an induced matching $M \in \mathcal{M}$ containing the vertex $v_i$ can only contain one vertex in $N_i$. Suppose this is not true, then $M$ contains distinct $x, y \in N_i$, without loss of generality we may assume that $v_i x$ and $yz$ are edges of $M$ for some vertex $z$. Since $v_i$ and $y$ are adjacent, this violates the assumption that $M$ is induced. It follows immediately from this observation that if $S$ is a subset of $\{v_1, \cdots, v_k\}$ and a matching $M \in \mathcal{M}$ contains precisely those vertices in $S$, then $M$ has at least $r-|S|$ edges inside the set $V \setminus (\bigcup_{i \in S} N_i)$. Now we let $M$ be a matching chosen uniformly at random from $\mathcal{M}$, and let $S$ be the subset of $\{v_1, \ldots, v_k\}$ that $M$ contains.
Define $Y_M$ to be a random subset $Y_M=V \setminus (\bigcup_{i \in S} N_i)$. Observe that the probability that a random matching in $\mathcal{M}$ contains a vertex $v_i$ is equal to $d_i/t$. Therefore, using $d_i \ge n_i$, 
\begin{align}\label{eqn_exp}
\mathbb{E} [|Y_M|] &= n - \sum_{i=1}^k \mathbb{P}(v_i \in M) \cdot |N_i| \le n -\sum_{i=1}^k n_i^2/t \nonumber\\
&\le n-1/t \cdot (e(G)-e(V_k)) \le n - 1/t \cdot (rt- n^2 / \log n)\\
&=(1-c+1/K) n\nonumber.
\end{align}

Now we partition the $t$ induced matchings from $\mathcal{M}$ into $2^k \le 2^{(\log n)/2} = n^{1/2}$ classes according to the set $V(M) \cap \{v_1, \cdots, v_k\}$. For a fixed constant $\alpha \in (0, 1)$, note that at most $\alpha t$ matchings in $\mathcal{M}$ come from those classes of sizes less than $\alpha Kn^{1/2}/(2 \log n)$, since otherwise bounding the number of matchings gives $\alpha t \le \alpha Kn^{1/2}/(2 \log n) \cdot n^{1/2}$, which results in a contradiction. Therefore the probability that a random matching $M$ belongs to one of the small classes $\leq \alpha$.
On the other hand, by Markov's inequality and \eqref{eqn_exp},
$$\mathbb{P}(|Y_M| \ge (1-c+1/K+\alpha)n) \le \frac{\mathbb{E}[|Y_M|]}{(1-c+1/K+\alpha)n} \le 1- \frac{\alpha}{1-c+1/K+\alpha}.$$
Since $c>1/5$ and $K>100/\varepsilon>100$, we have 
$$\frac{\alpha}{1-c+1/K+\alpha} > \frac{\alpha}{1-1/5+1/100+\alpha} > \alpha$$
if $\alpha \le 1/10$.
Therefore, by the union bound, we can find a class containing $t' \ge \alpha Kn^{1/2}/(2 \log n) = n^{1/2-o(1)}$ matchings, each of size at least $r'=r-(\log n)/2=(c-o(1))n$, on the same vertex set of size at most $n'=(1-c+1/K+\alpha)n$. When $c>1/5 + \varepsilon$, $K=100/\varepsilon$ and $\alpha$ is sufficiently small,
$$r'/n' = (c-o(1))/(1-c+1/K+\alpha) \ge 1/4.$$
However, from Theorem \ref{main_theorem}, we know that in this case $t'$ is at most $(6+o(1)) \log n' \le (6+o(1)) \log n$, which contradicts that $t' \ge n^{1/2-o(1)}$. \qed\\

\medskip
\noindent \textbf{Remark.} When the graph $G$ is regular or nearly regular, 
we can prove the $O(n/\log n)$ upper bound for a wider range $1/6+\varepsilon \le c < 1/4$. It is easy to see that every vertex has degree about $(2rt)/n=2ct$ in this case. We start with all the vertices uncovered and in each step we pick a vertex covering the most number of uncovered vertices. Let $X_i$ be the set of uncovered vertices in the $i$-th step. The sum of degrees of vertices in $X_i$ is equal to $2ct|X_i|$. Therefore, there exists a vertex $v_i$ covering at least $2ct|X_i|/n$ vertices of $X_i$. After $k$ steps, there are at most $n(1-2ct/n)^k$ uncovered vertices left. When $k=(\log n)/2$, $t>Kn/\log n$ and $K$ is sufficiently large, this gives $o(n)$ uncovered vertices. Note that from the proof of Theorem \ref{theorem_below_n/4}, the probability that a random matching from $\mathcal{M}$ covers a fixed vertex is equal to $d_i/t$, which is about $2ct/t=2c$. Therefore the expected number of vertices that are uncovered by vertices of $\{v_1, \cdots, v_k\} \cap V(M)$ is equal to $n-2c(n-o(n))=(1-2c+o(1))n$. Similarly as before we find a large class of induced matchings which are of size $(c-o(1))n$ and on the same vertex set of size at least $(1-2c+o(1))n$. When $c>1/6 + \varepsilon$, we are able to apply Theorem \ref{main_theorem} and derive a contradiction.\\

\noindent \textbf{Proof of Theorem \ref{theorem_just_below_n/4}:} 
Let $b=10^{-9}$. Suppose for contradiction $t>2n/(C \log n)$ and $n$ is sufficiently large. Here $C=C(n)$ is a function that tends to infinity slowly when $n \rightarrow \infty$, which we will choose later. By Lemma \ref{lemma_bipartite}, $F=G \times K_2$ is a $(2r, t)$-RS graph with parts $A$ and $B$, each of order $n$. Denote by $\mathcal{M}$ the collection of $t$ induced matchings that make up $F$. The number of edges of $F$ is $2rt=(1/2-2b)tn$. We will pick a sequence of disjoint subsets $A_1, \cdots, A_s$ from $A$, such that for each $i$, $|A_i|=Ct$, and $s$ is the smaller integer so that  $sCt \ge n/2$
(thus $s \sim n/(2Ct) < (\log n)/4$), moreover there is a subcollection of induced matchings $\mathcal{M}_i \subset \mathcal{M}$ such that $|\mathcal{M}_i| \ge t/16$ and each matching in $\mathcal{M}_i$ contains at most $qCt$ vertices from $A_i$, where $q=1/100$.

We first show that such a sequence of disjoint subsets would complete the proof. So suppose such a sequence of disjoints subsets exists. Define an auxiliary bipartite graph $X$ with first part being $[s]=\{1, \cdots, s\}$, and the second part being the matchings in $\mathcal{M}$, where $(i, M)$ is an edge of $X$ if and only if $M \in \mathcal{M}_i$. The minimum degree of vertices in $[s]$ is at least $t/16$, therefore the average degree of vertices in $\mathcal{M}$ is at least $(t/16)s/|\mathcal{M}|=s/16$. We can find a complete bipartite subgraph in $X$ with one part $|S|=s/16$, and the other part $\mathcal{M}' \subset \mathcal{M}$ with 
$$|\mathcal{M}'| \ge |\mathcal{M}|/\binom{s}{s/16} \ge |\mathcal{M}|/2^s \ge t/2^{(\log n)/4}>n^{1/2}.$$
Such a subgraph exists because the number of stars with a center in some $M \in \mathcal{M}$ and $s/16$ leaves in $[s]$ is equal to 
$\sum_{M \in \mathcal{M}} \binom{d_M}{s/16}$, where $d_M$ is the degree of $M$. By convexity and $\sum_M d_M \ge |\mathcal{M}|(s/16)$, this number is at least $|\mathcal{M}|$. Therefore, there exists such a subset $S$ of $[s]$ of size $s/16$, whose vertices have at least $|\mathcal{M}|/\binom{s}{s/16}$ common neighbors in $\mathcal{M}$.

Let $Y$ be the union of $A_i$ with $i \in S$, thus $|Y|=Ct(s/16) \ge n/32$. Each matching in $\mathcal{M}'$ contains at most $qCts/16 \le n/3200$ elements of $Y$; thus it contains at least $2r-n/3200=(1/2-2b-1/3200)n > 0.499n$ edges going between $A \setminus Y$ (which has size at most $31n/32$) and $B$ (which has size $n$). Note that $0.499n> \frac{1}{4}|B \cup (A \setminus Y)|=\frac{63}{128}n$. Hence, by Theorem \ref{double_count} there can only be a constant number of such induced matchings, a contradiction. So it suffices to show the existence of the desired disjoint subsets $A_1, \cdots, A_s$ of $A$.

We first show that the graph $F$ is nearly regular. Let $p=1/2000$, and $A=A' \cup A''$ be the partition of $A$ where $v \in A'$ if and only if $|d_F(v)-(t+1)/2| \le p(t+1)$. We will show that $A''$ is quite small. Consider the paths of length two with both end points in $A$. The number of such paths is equal to 
\begin{align}\label{2-path}
\sum_{v \in B} \binom{d_F(v)}{2} &=\frac{1}{2} \sum_{v \in B}d_F(v)^2 - rt \ge \frac{1}{2} \sum_{v \in B} (2rt/n)^2 -rt \nonumber\\
&=2(\frac{1}{4}-b)^2 t^2 n - (\frac{1}{4}-b)tn > (\frac{1}{8}-2b)(t+1)^2 n.
\end{align} 
Note that since $F$ is a $(2r, t)$-RS graph, every edge of $F$ has sum of degrees of its two vertices at most $t+1$. 
In particular, for any $v \in A$ and every its neighbor $u \in B$ we have $d(u) \leq (t+1-d(v))$.
Therefore, every vertex $v \in A$ is the end point of at most $d(v)(t+1-d(v)) \leq (t+1)^2/4$ paths of length two in $F$. Moreover, each vertex in $A''$ is the end point of at most $(\frac{1}{4}-p^2)(t+1)^2$ paths of length $2$ in $F$. Therefore the number of paths of length two with both endpoints in $A$ is at most $\frac{1}{2}(n-4p^2|A''|)(t+1)^2/4$. This estimate together with \eqref{2-path} implies that $|A''|<(4b/p^2)n = 2n/125$.

Similarly, let $B=B' \cup B''$ be the partition of $B$ where $B'$ consists of vertices which have degree within $p(t+1)$ of $(t+1)/2$, analogously we get $|B''|<2n/125$. Let $A^*$ be those vertices in $A$ with at least $t/4$ neighbors in $B''$.  Since the number of edges coming out of $B''$ is at most $|B''|t$, then $A^*$ has at most 
$|B''|t/(t/4)=4|B''|<8n/125$ elements. Thus there are at least $n-2n/125-8n/125 \ge 0.9n$ elements of $A$ neither in $A''$ nor $A^*$. Let $a_1$ be one such element in $A$, so $a_1$ is adjacent to at most $t/4$ neighbors in $B''$. Let $B_1$ be those vertices in $B'$ adjacent to $a_1$. So $|B_1| \ge (1/2-p)(t+1)-t/4 > t/8$, and $|B_1| \le t$. Let $A^{\#}$ be those vertices in $A$ adjacent to at least one vertex in $B_1$. The number of edges coming out from $B_1$ is at least $|B_1|(1/2-p)(t+1) \ge (t/8)(1/2-p)(t+1) > t^2/20$. If $|A^{\#}|<Ct$, then we would have a bipartite graph between $A^{\#}$ and $B_1$ with parts of size at most $Ct$ and $t$ respectively but with at least $t^2/20$ edges. This graph is an RS-graph on at most $(C(n)+1)t$ vertices with at least $t^2/20$ edges which can be decomposed into $t$ induced matchings. It is well known that this leads to a graph $R$ on $\sim C(n)t$ vertices with $\Omega(t^2)$ edges, in which every edge is contained in exactly one triangle (see e.g. Proposition 3.1 in \cite{alon-moitra-sudakov}). We take $C(n)$ to be  $2^{\Omega(\log^{*} n)}$, where $\log ^{*}$ is the iterated logarithm. The best known bound on the triangle removal lemma \cite{fox} (see also \cite{cofo} for an alternative proof) guarantees that such a graph $R$ cannot exist. We therefore have $|A^{\#}| \ge Ct$. Let $A_1$ be an arbitrary subset of $A^{\#}$ of size $Ct$. 

Let $\mathcal{M}'$ be the collection of matchings in $\mathcal{M}$ containing $a_1$, and $\mathcal{M}''$ be the matchings not containing $a_1$, then $|\mathcal{M}'| \ge (1/2-p)(t+1)$ and $|\mathcal{M}''| \le t-(1/2-p)(t+1)$. Each vertex $v \in B_1$ is in at least $(1/2-p)(t+1)$ matchings, and all but one of them is in $\mathcal{M}''$. Thus a vertex $u$ in $A_1$ (which is adjacent to some $v \in B_1$) cannot be contained in these matchings except for one, so $u$ is in at most $t-(1/2-p)(t+1)-(1/2-p)(t+1)+1=2p(t+1)$ matchings in $\mathcal{M}''$. Therefore the number of matchings in $\mathcal{M}''$ that contain at least $qCt$ elements of $A_1$ is at most $2p(t+1)(Ct)/(qCt)=2p(t+1)/q<t/8 \le |\mathcal{M}''|/2$. Therefore at least half the matchings in $\mathcal{M}''$ contain at most $qCt$ elements of $A_1$. Let $\mathcal{M}_1$ be the set of these matchings, we have $|\mathcal{M}_1| \ge t/16$.

We have thus pulled out the desired sets $A_1$ and $\mathcal{M}_1$ and finished the first step. In the second step, we find $A_2$ in $A \setminus A_1$ and $\mathcal{M}_2$ in a similar fashion. After the $i$-th step, we have pulled out disjoint sets $A_1, \cdots , A_i$ and also found sets of matchings $\mathcal{M}_1, \cdots ,\mathcal{M}_i$.
Note that $|A'' \cup A^* \cup A_1 \cup \cdots \cup A_i| \le 2n/125+8n/125+sCt \le 0.7n$. Let $A_0=A \setminus (A'' \cup A^* \cup A_1 \cup \cdots \cup A_i)$, then $|A_0| \ge 0.3n$. Note that by similar arguments, every vertex in $A_0$ is adjacent to at least $t/8$ neighbors in $B'$. Moreover let $e_v$ be the number of paths $vuv'$ with $v, v' \in A_0$ and $u \in B'$. Note that $e_v, v \in A_0$ counts precisely the number of edges between the neighbors of $v$ in $B'$ and the second neighborhood of $v$ in $A_0$.
For $u \in B'$, let $d'_u$ be the number of vertices in $A_0$ adjacent to $u$, then
\begin{align*}
\sum_v e_v = 2 \sum_{u \in B'} \binom{d'_u}{2} \ge 2 \binom{e(A_0, B')/|B'|}{2} |B'|
\end{align*}
Therefore there exists a vertex $a_{i+1} \in A_0$, such that 
$$e_{a_{i+1}} \ge 2 \binom{e(A_0, B')/|B'|}{2} |B'|/|A_0| = \frac{e(A_0, B')}{|A_0|} \left(\frac{e(A_0, B')}{|B'|}-1 \right).$$

Since every vertex in $A_0$ is adjacent to at least $t/8$ vertices in $B'$, we have $e(A_0, B')/|A_0| \ge t/8$, and $e(A_0, B') \ge 0.3nt/8$. Hence $e_{a_{i+1}} \ge (t/8)(0.3nt/8n-1) \ge t^2/300$. Once again the triangle removal lemma implies that
we can pick $A_{i+1}$ to be any $Ct$ vertices from $A_0$ that are second neighbors of $a_{i+1}$. Let $\mathcal{M}_{i+1}$ be the collections of matchings in $\mathcal{M}$ not containing $a_{i+1}$ and containing at most $qCt$ elements of $A_{i+1}$, similarly we have $|\mathcal{M}_{i+1}| \ge t/16$. In the end we obtain the disjoint subsets $A_1, \cdots, A_s$ and collections of matchings $\mathcal{M}_1, \cdots, \mathcal{M}_s$ as desired. 
\qed\\

\section{Concluding Remarks} \label{section_concluding}
In this paper we consider the problem of determining the maximum number of induced matchings $t$ as a function in the size of matchings $r$ and number of vertices $n$, and settle the problem for $r=cn$ when $c>1/4$ and $c=1/4$. Several intriguing problems remain open.
\begin{list}{\labelitemi}{\leftmargin=1em}
\item When $r=n/4$, we prove that $t \le (6+o(1)) \log_2 n$, while the hypercube construction and its refinement only gives a $t= (2+o(1)) \log_2 n$ lower bound. For regular graphs, we know this lower bound is tight.
It seems plausible that for $r, n$ satisfying $r=n/4$ and $n$ tending to infinity, we always have $t \le 2 \log_2 n +2$.
The expansion part of our proof of Theorem \ref{main_theorem} is pretty robust in that it uses nothing about the number $t$ of matchings. However, there are two places we lost a constant factor in the proof. We first pass from $G$ to $H$ (resp. $H'$), keeping only edges whose sum of degrees is at least $t$ (resp. $t+1$). Maybe it is possible to somehow incorporate negative edges (i.e., whose sum of degrees of its two vertices is less than $t$) in the argument. In the second step, we pass from $H$ (or $H'$) to its subgraph by keeping only high degree vertices. It would be great to do a more careful analysis and eventually improve the constant before $\log_2n$ from $6$ to $2$.

\item When the size of matchings is close to $n/4$, there are two quite different constructions of Ruzsa-Szemer\'edi graphs. One is
the Kneser graph $KG(2k+1, k)$ with $k \sim \frac{1}{2} \log_2 n$, which is an $(n/4 + \Theta(n/\log n), (1+o(1))\log_2 n)$-RS graph. The other is the hypercube $\{0,1\}^{\log_2 n}$, which is an $(n/4, 2 \log_2 n)$-RS graph. Can we find a family of $(r, t)$-RS graphs that bridge between these two examples, say with $r=c \log_2 n$ for some $c \in [1, 2]$, and $t-n/4 = \Omega (\log n)$?

\item In Section \ref{section_r<n/4} we obtained an upper bound $O(n/ \log n)$ on the number of matchings when $c$ is between $1/5$ and $1/4$, and improved the bound to $n/((\log n)  2^{\Omega(\log^* n)})$ for $c>1/4-b$ for some constant $b>0$. However the best construction \cite{fischer-et-al} so far only gives $n^{\Omega(1/ \log\log n)}$. It remains a challenging question to close this gap for this range, and in general for all $r$ linear in $n$.

\item Another interesting open problem is to maximize the size of induced matchings, when the number of matchings is linear in $n$. It is somewhat complementary to the questions we studied in this paper. The Ruzsa-Szemer\'edi theorem
\cite{ruzsa-szemeredi} gives an upper bound $o(n)$, which in turn gives a similar bound for Roth's theorem. On the other hand, 
for Roth's theorem, much better results are known. The current best upper bound on the size of the largest subset of $[n]$
with no $3$-term arithmetic progression is of the form $O(n/\log^{1-o(1)} n)$, first proved by Sanders \cite{sanders} using Fourier analysis on Bohr sets and subsequently improved slightly further by Bloom \cite{bloom}. It would be very interesting to see whether a similar bound holds for the Ruzsa-Szemer\'edi theorem. Moreover, the best known lower bound for both problems is
still essentially $n/e^{O(\sqrt{\log n})}$ obtained by Behrend \cite{behrend} about seventy years ago. Closing the gap between these bounds
is an important open problem.
\end{list}

\end{document}